\pdfoutput=1
\RequirePackage{ifpdf}
\ifpdf 
\documentclass[pdftex]{sigma}
\else
\documentclass{sigma}
\fi

\DeclareMathOperator{\Ker}{Ker}
\DeclareMathOperator{\Image}{Im}
\DeclareMathOperator{\rank}{rank}

\newtheorem*{lemma-nonnumbered}{Lemma}
\newtheorem*{xpr}{Experimental result}
\theoremstyle{definition}
\newtheorem{cnv}{Convention}

\begin{document}


\renewcommand{\PaperNumber}{053}

\FirstPageHeading

\ShortArticleName{Parameterizing the Simplest Grassmann--Gaussian Relations for Pachner Move 3--3}

\ArticleName{Parameterizing the Simplest Grassmann--Gaussian\\ Relations for Pachner Move 3--3}

\Author{Igor G.~KOREPANOV and Nurlan M.~SADYKOV}

\AuthorNameForHeading{I.G.~Korepanov and N.M.~Sadykov}

\Address{Moscow State University of Instrument Engineering and Computer Sciences,\\ 20 Stromynka Str., Moscow
107996, Russia} \Email{\href{mailto:paloff@ya.ru}{paloff@ya.ru},
\href{mailto:nik56-74@mail.ru}{nik56-74@mail.ru}}

\ArticleDates{Received May 15, 2013, in f\/inal form August 08, 2013; Published online August 13, 2013}

\Abstract{We consider relations in Grassmann algebra corresponding to the four-dimen\-sional Pachner move
3--3, assuming that there is just one Grassmann variable on each 3-face, and a~4-simplex weight is
a~Grassmann--Gaussian exponent depending on these variables on its f\/ive 3-faces.
We show that there exists a~large family of such relations; the problem is in f\/inding their
algebraic-topologically meaningful parameterization.
We solve this problem in part, providing two nicely parameterized subfamilies of such relations.
For the second of them, we further investigate the nature of some of its parameters: they turn out to
correspond to an exotic analogue of middle homologies.
In passing, we also provide the 2--4 Pachner move relation for this second case.}


\Classification{15A75; 57Q99; 57R56}

\section{Introduction}

Discrete topological f\/ield theories~-- specif\/ically, f\/ield theories on piecewise linear (PL)
manifolds~-- are def\/initely a~challenging research subject.
As there are now many interesting topological quantum f\/ield theories (TQFT's) in three dimensions, it
looks reasonable to concentrate on the four-dimen\-sional case.
Such a~theory is expected to bring about interesting results both by itself and when compared with the
existing theories on smooth manifolds.

As explained, for instance, in~\cite[Section~1]{Lickorish}, it makes sense f\/irst to construct
algebraic relations corresponding to \emph{Pachner moves}.
And the simplest nontrivial relations of such kind arise, as we believe, in \emph{Grassmann algebras}.
In three dimensions, a~relation corresponding to Pachner move 2--3 is often called \emph{pentagon
relation}, and there are some Grassmann-algebraic constructions for pentagon relation, presented, in
particular, in paper~\cite{KS}.
As we hope to demonstrate here, the four-dimensional case has its own specif\/ic beauty; it is more
complicated but also yields to systematic investigation.

If we consider an \emph{ansatz}~-- a~(tentative) specif\/ic form of quantities or expressions entering in
our relations, and consider the relations as equations for the ansatz parameters, and if our ansatz is
simple enough, then it may happen that the \emph{existence} of many solutions for such equations follows
already from parameter counting.

In this paper, we take the simplest possible form~\eqref{33} of Grassmann-algebraic relation corresponding
to Pachner move 3--3~-- with just one Grassmann variable on each 3-face, and further assume that the
Grassmann weight of a~4-simplex has the form of a~Grassmann--Gaussian exponent, depending on the f\/ive
variables on the 3-faces.
A heuristic parameter count shows that there exists a~large~-- and intriguing~-- family of relations of
such form.
We prefer to go further and prove the rigorous Theorem~\ref{th:9c}, formulated in terms of \emph{isotropic
linear spaces of Grassmann differential operators} annihilating our Grassmann--Gaussian exponents.
In doing so, we not only prove the existence of the 4-simplex weights satisfying the 3--3 relations, but
discover some interesting operators (namely,~\eqref{g} and~\eqref{h}) that may deserve further
investigation; at least, they have an elegant form (namely,~\eqref{g1} and~\eqref{h1}) in one specif\/ic
case.

Having proved our Theorem~\ref{th:9c}, we are naturally led to the problem of f\/inding an
\emph{algebraic-topologically meaningful parameterization} of our Grassmann weights, which would enable us
to move further and construct topological f\/ield theories.
In the present paper, we make two steps in this direction by presenting two explicitly~-- and nicely~--
parameterized \emph{subfamilies} of such weights, found largely by guess-and-try method.
The f\/irst subfamily resembles the (more cumbersome) constructions in~\cite{11S,exo}~-- both are related
to \emph{exotic homologies}.
The striking new fact is, however, that this is now only a~subfamily of something mysterious, on whose
nature only our parameterized second family sheds some additional light.

Some of the results of this paper have f\/irst appeared, in a~preliminary form, in the
preprint~\cite{2-cocycles}.

Below,
\begin{itemize}
\itemsep=0pt
\item
in Section~\ref{s:G}, we recall the basic def\/initions from the theory of Grassmann algebras and
Berezin integral,
\item
in Section~\ref{s:P}, we recall the four-dimensional Pachner moves, mainly moves 3--3 and 2--4 with
which we will be dealing in this paper, and introduce some notational conventions,
\item
in Section~\ref{s:33gen}, we introduce a~3--3 relation for Grassmann 4-simplex weights.
First, we do it in a~general form, then we specialize the weights to be Grassmann--Gaussian exponents and
explain their connection with isotropic spaces of Grassmann dif\/ferential operators,
\item
in Section~\ref{s:p}, based on these isotropic spaces, we show the existence of a~vast family of
4-simplex weights satisfying the 3--3 relation.
The way we do it is constructive; what lacks in it is a~parameterization for this whole family relevant for
algebraic-topological applications,
\item
in Section~\ref{s:33spe}, we present two subfamilies of Grassmann 4-simplex weights satisfying the
3--3 relation where such parameterization has been obtained,
\item
in Section~\ref{s:h}, we do some preparational work in order to expose some exotic-homological
structures lying behind the second of the mentioned subfamilies.
Namely, we introduce, for a~given triangulated four-manifold, a~sequence of two linear mappings~--
supposedly a~fragment of an exotic chain complex, prove their chain property (their composition vanishes),
and present computational evidence showing that they provide an exotic analogue of usual middle (i.e.,
second) homologies, and
\item
in Section~\ref{s:24spe}, guided by the fact that the mentioned exotic-homological structures
manifest themselves more clearly for the Pachner move 2--4, we present the relations corresponding to this
move, study a~new factor~-- edge weight~-- appearing in these relations, and then formulate the relations
for both moves 3--3 and 2--4 using these exotic-homological terms.
\end{itemize}

\section{Grassmann algebras and Berezin integral}
\label{s:G}

In this paper, a~\emph{Grassmann algebra} is an associative algebra over the f\/ield $\mathbb C$ of complex
numbers, with unity, generators $x_i$~-- also called Grassmann variables~-- and relations
\begin{gather*}
x_i x_j=-x_j x_i.
\end{gather*}
This implies that, in particular, $x_i^2 =0$, so each element of a~Grassmann algebra is a~polynomial of
degree $\le 1$ in each $x_i$.

The \emph{degree} of a~Grassmann monomial is its total degree in all Grassmann variables.
If an algebra element consists of monomials of only odd or only even degrees, it is called \emph{odd} or,
respectively, \emph{even}.
If all the monomials have degree 2, we call such element a~Grassmannian \emph{quadratic form}.

The \emph{exponent} is def\/ined by its usual Taylor series.
We call the exponent of a~quadratic form \emph{Grassmann--Gaussian exponent}.
Here is an example of it:
\begin{gather*}
\exp(x_1x_2+x_3x_4)=1+x_1x_2+x_3x_4+x_1x_2x_3x_4.
\end{gather*}

There are two kinds of derivations in a~Grassmann algebra: left derivative $\frac{\partial}{\partial x_i}$
and right derivative $\frac{\overleftarrow{\partial}}{\partial x_i}$, with respect to a~Grassmann variable
$x_i$.
These are $\mathbb C$-linear operations in Grassmann algebra def\/ined as follows.
Let $f$ be an element not containing variable $x_i$, then
\begin{gather}
\label{Gd1}
\dfrac{\partial}{\partial x_i}f=f\dfrac{\overleftarrow{\partial}}{\partial x_i}=0,
\end{gather}
and
\begin{gather}
\label{Gd2}
\dfrac{\partial}{\partial x_i}(x_if)=f,
\qquad
(fx_i)\dfrac{\overleftarrow{\partial}}{\partial x_i}=f.
\end{gather}

From~\eqref{Gd1} and~\eqref{Gd2}, the following \emph{Leibniz rules} follow: if $f$ is either even or odd,
then
\begin{gather}
\label{L}
\dfrac{\partial}{\partial x_i}(fg)=\dfrac{\partial}{\partial x_i}f\cdot g+\epsilon f\dfrac{\partial}
{\partial x_i}g,
\qquad
(gf)\dfrac{\overleftarrow{\partial}}{\partial x_i}=g\cdot f\dfrac{\overleftarrow{\partial}}{\partial x_i}
+\epsilon g\dfrac{\overleftarrow{\partial}}{\partial x_i}f,
\end{gather}
where $\epsilon=1$ for an even $f$ and $\epsilon=-1$ for an odd $f$.

The \emph{Grassmann--Berezin calculus} of anticommuting variables is in many respects parallel to the usual
calculus, see~\cite{B} and especially~\cite{B-super}.
Still, there are some peculiarities, and one of them is that the integral in a~Grassmann algebra is the
same operation as derivative; more specif\/ically, \emph{Berezin integral} in a~variable $x_i$ is
def\/ined, traditionally, as the \emph{right} derivative w.r.t.~$x_i$.
Independently, Berezin integral is def\/ined as follows: it is a~$\mathbb C$-linear operator in Grassmann
algebra satisfying
\begin{gather*}
\int\mathrm dx_i=0,
\qquad
\int x_i\,\mathrm dx_i=1,
\qquad
\int gh\,\mathrm dx_i=g\int h\,\mathrm dx_i,
\end{gather*}
where $g$ does not contain $x_i$.
Multiple integral is def\/ined according to the following Fubini rule:
\begin{gather}
\label{mB}
\idotsint f\,\mathrm dx_1\,\mathrm dx_2\,\cdots\,\mathrm dx_n=\int\left(\cdots\int\left(\int f\,\mathrm dx_1\right)\mathrm dx_2\cdots\right)\mathrm dx_n .
\end{gather}
In ``dif\/ferential'' notations, integral~\eqref{mB} is
\begin{gather}
\label{dmB}
f\dfrac{\overleftarrow{\partial}}{\partial x_1}\dfrac{\overleftarrow{\partial}}{\partial x_2}
\cdots\dfrac{\overleftarrow{\partial}}{\partial x_n}.
\end{gather}

\section{Pachner moves in four dimensions}
\label{s:P}

Pachner moves~\cite{Pachner} are elementary local rebuildings of a~manifold triangulation.
A triangulation of a~PL manifold can be transformed into any other triangulation using a~f\/inite sequence
of Pachner moves.

In four dimensions, each Pachner move replaces a~cluster of 4-simplices with a~cluster of some other
4-simplices, occupying the same place in the triangulation and having the same boundary.
There are f\/ive (types of) Pachner moves in four dimensions: $3\to 3$, $2\leftrightarrow 4$, and
$1\leftrightarrow 5$, where the numbers indicate how many 4-simplices have been withdrawn and how many have
replaced them.
As the withdrawn and the replacing clusters of 4-simplices have the same common boundary, we can glue them
together in a~natural way (forgetting for a~moment about the rest of the manifold); then, for all Pachner
moves, they must form together a~sphere $S^4$ triangulated in six 4-simplices as the boundary of
a~5-simplex, which we denote $\partial \Delta^5$.
More details and a~pedagogical introduction can be found in~\cite{Lickorish}.

More traditional notations for the mentioned moves are 3--3, 2--4, 4--2, 1--5, and 5--1; we will be using
these notations as well.

Move 3--3 is, in some informal sense, central: experience shows that if we have managed to f\/ind an
algebraic formula whose structure can be regarded as ref\/lecting the structure of the move, then we can
also f\/ind (usually more complicated) formulas corresponding to the other Pachner moves.
This may be compared to the three-dimensional case, where the popular ``pentagon relation'' often
corresponds to the ``central'' three-dimensional Pachner move 2--3, while it is believed that, having done
something interesting with this pentagon equation, one will be also able to work with the move 1--4.

We call the initial cluster of 4-simplices in a~move the \emph{left-hand side} (l.h.s.) of that move, and
the f\/inal cluster~-- its \emph{right-hand side} (r.h.s.).
All moves in this paper will involve six vertices denoted $i=1,\dots,6$.
Below are some more details.

\subsection[Move $3\to 3$]{Move $\boldsymbol{3\to 3}$}

It transforms, in the notations used in this paper, the cluster of three 4-simplices 12345, 12346 and 12356
situated around the 2-face 123 into the cluster of three other 4-simplices, 12456, 13456 and 23456,
situated around the 2-face 456.
The inner 3-faces (tetrahedra) are 1234, 1235 and 1236 in the l.h.s., and 1456, 2456 and 3456 in the r.h.s.
The boundary of both sides consists of nine tetrahedra listed below in table~\eqref{table}.

There are no inner edges (1-faces) or vertices (0-faces) in either side of this move.

\subsection[Moves $2\leftrightarrow 4$]{Moves $\boldsymbol{2\leftrightarrow 4}$}

We describe move $2\to 4$; move $4\to 2$ is its inverse.
Move $2\to 4$ replaces, in the notations of this paper, the cluster of two 4-simplices 12345 and 12346 with
the cluster of four 4-simplices 12356, 12456, 13456 and 23456.
The boundary of both sides consists of eight tetrahedra 1235, 1236, 1245, 1246, 1345, 1346, 2345 and 2346.

In the l.h.s., there is one inner tetrahedron 1234, no inner 2-faces and no inner edges.

In the r.h.s., there are six inner tetrahedra 1256, 1356, 1456, 2356, 2456 and 3456, and four inner
2-faces: 156, 256, 356 and 456.
It turns out especially important~-- see Subsection~\ref{ss:w}~-- that there \emph{is} one inner edge,
namely 56, in the r.h.s.

There are no inner vertices in either side of this move.

\subsection[Moves $1\leftrightarrow 5$]{Moves $\boldsymbol{1\leftrightarrow 5}$}

We don't work with these moves in this paper, so we only indicate that the move $1\to 5$ adds a~new vertex
$6$ inside the 4-simplex $12345$, thus dividing it into f\/ive new 4-simplices.
Move $5\to 1$ is, of course, its inverse.

\subsection{A few conventions}

Any side of a~Pachner move, as well as a~single 4-simplex, is a~triangulated four-manifold with boundary.
In Section~\ref{s:h}, we also consider an arbitrary \emph{orientable} triangulated four-manifold.
Here are some our conventions concerning manifolds, their simplices, and also some complex parameters
appearing in our theory, like \emph{vertex coordinates} (see Section~\ref{s:33spe}).
\begin{cnv}
All manifolds in this paper are assumed to be \emph{oriented}.
In the case of Pachner moves, the orientation is def\/ined so that, for the 4-simplex 12345, it is given by
this order of its vertices.
\end{cnv}
\begin{cnv}
\label{cnv:num}
We denote by $N_k$ the number of $k$-simplices in a~triangulation, and by $N'_k$ the number of \emph{inner}
$k$-simplices.
Vertices are numbered from 1 through $N_0$ (as we have already done for Pachner moves, where $N_0=6$).
\end{cnv}
\begin{cnv}
When simplices are written in terms of their vertices, these go in the increasing order of their numbers
(again, as we have already done).
\end{cnv}
\begin{cnv}
\label{cnv:set}
If the order of vertices is unknown, we use the following notation.
Let there be, for instance, a~4-simplex whose set of vertices is $\{i,j,k,l,m\}$, then we denote it,
omitting the commas for brevity, as $\{ijklm\}$.
\end{cnv}
\begin{cnv}
The complex parameters appearing in our theory~-- to be exact, the eighteen parameters in
Section~\ref{s:p} and vertex coordinates introduced in Section~\ref{s:33spe}~-- lie in the
\emph{general position} with respect to any considered algebraic formula, unless the opposite is explicitly
stated.
For instance, there is no division by zero in formula~\eqref{tg}.
Moreover, concerning vertex coordinates, all functions of them in this paper are \emph{rational}, so the
reader can assume that the coordinates are \emph{indeterminates} over $\mathbb C$ (and we extended $\mathbb
C$ to the relevant f\/ield of rational functions).
\end{cnv}

\section[Relation 3-3 with Grassmann-Gaussian exponents: generalities]{Relation 3--3 with Grassmann--Gaussian exponents:\\ generalities}
\label{s:33gen}

\subsection{The form of relation 3--3}

The Grassmann-algebraic Pachner move relations for move 3--3, considered in this paper, have the following
general form:
\begin{gather}
\iiint\mathcal W_{12345}\mathcal W_{12346}\mathcal W_{12356}\,\mathrm dx_{1234}\,\mathrm dx_{1235}
\,\mathrm dx_{1236}
\nonumber
\\
\qquad
{}=\mathrm{const}\iiint\mathcal W_{12456}\mathcal W_{13456}\mathcal W_{23456}\,\mathrm dx_{1456}
\,\mathrm dx_{2456}\,\mathrm dx_{3456}.
\label{33}
\end{gather}
Here Grassmann variables $x_{ijkl}$ are attached to all 3-faces, i.e., tetrahedra $t=ijkl$; the
\emph{Grassmann weight} $\mathcal W_{ijklm}$ of a~4-simplex $u=ijklm$ depends on (i.e., contains) the
variables on its 3-faces, e.g., $\mathcal W_{12345}$ depends on $x_{1234}$, $x_{1235}$, $x_{1245}$,
$x_{1345}$ and $x_{2345}$.
Also, $\mathcal W_u$ may depend on parameters attached to the 4-simplex $u$ or/and its subsimplices.
The integration goes in variables on \emph{inner} three-faces in the corresponding side of Pachner move,
while the result depends on the variables on boundary faces.
Finally, $\mathrm{const}$ in the right-hand side is a~numeric factor.

Formula~\eqref{33} appears to give the simplest possible form for a~Grassmann-algebraic relation imitating
the 3--3 move.

\begin{remark}
We expect formula~\eqref{33} to be the f\/irst step in constructing a~fermionic topological f\/ield theory,
based on our previous experience.
As yet, the most detailed construction of such kind can be found in paper~\cite{0806}, although for
a~three-dimensional case and a~simpler theory.
Note especially
\begin{itemize}\itemsep=0pt
\item
the analogy between either side of~\eqref{33} and~\cite[formula~(10)]{0806} describing what happens when we
glue together manifolds with boundary,
\item
also the analogy between either side of~\eqref{33} if the $\mathcal W_u$ are as in formula~\eqref{gg}
below, and~\cite[formula~(8)]{0806} for a~multi-component fermionic state-sum invariant of a~manifold with
boundary,
\item
and the discussion in~\cite[Section~5]{0806} of the modif\/ication of Atiyah's TQFT axioms for
anti-commuting variables.
\end{itemize}
\end{remark}

\begin{remark}
As an example of a~more heavyweight relation corresponding to the same Pachner move 3--3, we can
cite~\cite[formula (38)]{11S}, where, in particular, \emph{two} Grassmann variables live on each
tetrahedron.
\end{remark}

Further simplif\/ication is achieved by using \emph{Grassmann--Gaussian exponents} (which corresponds to
\emph{free fermions} in physical language) and assuming that
\begin{gather}
\label{gg}
\mathcal W_u=\exp\left(\frac{1}{2}\sum_{t_1,t_2\subset u}\alpha_{t_1t_2}^{(u)}x_{t_1}x_{t_2}\right),
\end{gather}
where $t_1$ and $t_2$ are two 3-faces of $u$, and $\alpha_{t_1t_2}^{(u)}\in\mathbb C$ are numeric
coef\/f\/icients, with the antisymmetry condition
\begin{gather}
\label{as}
\alpha_{t_1t_2}^{(u)}=-\alpha_{t_2t_1}^{(u)}.
\end{gather}
We hope to demonstrate in this paper that the relation~\eqref{33}
is interesting already in the case of
such exponents.

\subsection{Isotropic subspaces of operators}

The exponent~\eqref{gg} is characterized, up to a~factor that does not depend on those $x_t$ that enter in
it, by the equations
\begin{gather}\label{dandx}
\left(\partial_t-\sum_{t'\subset u}\alpha_{tt'}^{(u)}x_{t'}\right)\mathcal W_u=0
\qquad
\text{for all}
\quad
t\subset u,
\end{gather}
where we denote $\partial_t = \partial / \partial x_t$.
Generalizing the operators in the big parentheses in~\eqref{dandx}, we consider $\mathbb C$-linear
combinations of operators of left dif\/ferentiations and multiplying by Grassmann generators:
\begin{gather}
\label{bg}
d=\sum_t(\beta_t\partial_t+\gamma_tx_t),
\end{gather}
where $t$ runs over all 3-faces in a~given triangulated manifold.
\begin{remark}
\label{r:PWZ}
We are going to apply dif\/ferential equations like~\eqref{dandx} to analysing the integral
identity~\eqref{33}.
Such kind of analysis is a~well known and widely used tool, especially in the bosonic context.
For example, holonomic and $q$-holonomic functions are studied via the dif\/ference/dif\/fe\-ren\-tial
equations they satisfy, see~\cite{Bjork,Cartier,PWZ}.
\end{remark}

We regard the anticommutator of two operators~\eqref{bg} (def\/ined as $[A,B]_+=AB+BA$ for operators $A$
and $B$) as their \emph{scalar product}:
\begin{gather}
\label{sc}
\big\langle d^{(1)},d^{(2)}\big\rangle\stackrel{\rm def}{=}\big[d^{(1)},d^{(2)}\big]_+
=\sum_t\Big(\beta_t^{(1)}\gamma_t^{(2)}+\beta_t^{(2)}\gamma_t^{(1)}\Big).
\end{gather}
With this scalar product, operators~\eqref{bg} form a~\emph{complex Euclidean space}, while all polynomials
of these operators form a~\emph{Clifford algebra}.

Recall that an \emph{isotropic}, or \emph{totally singular}, subspace of a~complex Euclidean space $\mathbb
C^{2n}$ is such linear subspace where the scalar product identically vanishes.
We will need some basic facts about isotropic subspaces; for the reader's convenience, we formulate them as
the following Theorem~\ref{th:iso} and give it a~simple proof.
Much more interesting facts about Clif\/ford algebras and isotropic subspaces in Euclidean spaces can be
found, e.g., in~\cite{C}.

\begin{theorem}\label{th:iso}
Maximal isotropic spaces in complex Euclidean space $\mathbb C^{2n}$ have dimension~$n$.
The manifold of these maximal isotropic spaces~-- isotropic Grassmannian~-- splits up in two
connected components.
\end{theorem}

\begin{proof}
The f\/irst statement is an easy exercise.
The second can be proved as follows.
Let $V\subset \mathbb C^{2n}$ be a~maximal isotropic subspace.
For a~\emph{generic} orthonormal basis $\mathsf e_1,\dots,\mathsf e_{2n}$, the orthogonal projection of $V$
onto the space $W$ spanned by the f\/irst half of basis vectors, i.e., $\mathsf e_1,\dots,\mathsf e_n$,
coincides with the whole $W$.
Also, considering the manifold $B$ of orthonormal bases, it has two connected components (the determinant
of transition matrix is $1$ within a~component, and $-1$ between the components), and the same components
remain in $B\setminus S$~-- the result of taking away the set $S$ of non-generic (in the sense indicated
above) bases: the components cannot split any further, because $S$ has complex codimension $\ge 1$ and thus
real codimension $\ge 2$.
Taking some liberty, we call, in this proof, the components of $B$ \emph{orientations} of $\mathbb C^{2n}$.

We prefer to arrange basis vectors in a~column, and vector coordinates in a~row; so, an arbitrary vector in
$\mathbb C^{2n}$ is written like $
\begin{pmatrix}
\alpha_1 & \dots & \alpha_n
\end{pmatrix}
\begin{pmatrix}
\mathsf e_1
\\
\vdots
\\
\mathsf e_{2n}
\end{pmatrix}
$, and vectors in $\mathbb C^{2n}$ are identif\/ied with row vectors if a~basis is given.
For a~space $V$ and basis $\mathsf e_1,\dots,\mathsf e_{2n}$ such as in the previous paragraph, we can
represent $V$ as the linear span of the rows of the following matrix:
\begin{gather}\label{1o}
\begin{pmatrix}
\mathbf1_n&\mathrm i O
\end{pmatrix}
,
\end{gather}
where $\mathbf 1_n$ is the identity matrix and $O$ is an orthogonal matrix, both of sizes $n\times n$, and
$\mathrm i=\sqrt{-1}$.
The determinant of $O$ is either $1$ or $-1$, and its sign obviously cannot change within one component of
$B\setminus S$.
So, given a~f\/ixed orientation of $\mathbb C^{2n}$ and a~maximal isotropic space $V$, we get either $1$ or
$-1$ as $\det O$.

It remains to prove that there is no further splitting between maximal isotropic spaces.
Consider two such spaces, $V_1$ and $V_2$.
There exists a~generic, in the above sense, basis $\mathsf e_1,\dots,\mathsf e_{2n}$ for \emph{both} of
them.
Using this basis, $V_1$ and $V_2$ can be written in terms of matrices~\eqref{1o}, and the corresponding
orthogonal matrices $O_1$ and $O_2$ belong to the same connected component in the space of orthogonal
matrices.
\end{proof}

\begin{theorem}\label{th:5}\qquad

\begin{enumerate}\itemsep=0pt
\item
\label{th:5:i}
For a~given weight $\mathcal W_u$ of the form~\eqref{gg}, the operators $d$ satisfying equation
\begin{gather}
\label{dWu}
d\mathcal W_u=0
\end{gather}
form a~five-dimensional isotropic linear space.
\item
\label{th:5:ii}
For the set of equations~\eqref{dWu} corresponding to a~five-dimensional isotropic space $V$ of
operators~\eqref{bg} with $t$ running over the five $3$-faces of a~$4$-simplex $u$,
there exists a~nonzero Grassmann algebra element $\mathcal W_u$, containing only the Grassmann generators
$x_t$ and sa\-tis\-fying all these equations.
This $\mathcal W_u$ is determined by these equations uniquely up to a~numeric factor.
\item
\label{th:5:iii}
The element $\mathcal W_u$ from item~$\ref{th:5:ii}$ is even for one connected component of the set of
five-dimen\-sional isotropic spaces~$V$, and odd for the other.
In the first case, it is, for a~gene\-ric~$V$, a~Grassmann--Gaussian exponent~\eqref{gg}.
\end{enumerate}
\end{theorem}

\begin{proof}
 (i)~Five such linearly independent equations are already written in~\eqref{dandx}.
It follows from the antisymmetry~\eqref{as} of coef\/f\/icients~$\alpha_{tt'}$ that any two operators~$d$
written in the big parentheses in~\eqref{dandx} anticommute (including the case where they coincide).
This means that the scalar product~\eqref{sc} vanishes.

(ii)~We denote the ten-dimensional Euclidean space of all operators~\eqref{bg}, where $t\subset
u$, simply as $\mathbb C^{10}$.
There exists an orthogonal transform $O$ of $\mathbb C^{10}$ sending $V$ into the subspace generated by the
f\/ive $\frac{\partial}{\partial x_t}$, and to $O$ corresponds, according to the general theory, a~$\mathbb
C$-linear automorphism $B$ of the Grassmann algebra such that
\begin{gather*}
O\mathsf y=B\mathsf y B^{-1}
\qquad
\text{for}
\quad
\mathsf y\in\mathbb C^{10}.
\end{gather*}
As $\frac{\partial}{\partial x_t}1=0$ for all f\/ive~$t$, it follows that $\mathcal W_u = B^{-1} 1$ is
annihilated by all operators in $V$.
On the other hand, if there were two linearly independent~$\mathcal W_u$ annihilated by all operators in
$V$, it would follow that the two corresponding algebra elements $B \mathcal W_u$ would be annihilated by
all f\/ive~$\frac{\partial}{\partial x_t}$, but this holds only for the one-dimensional space of constants.

(iii)~A Zariski open set of even elements $\mathcal W_u$ is already provided, and it consists of
Grassmann--Gaussian exponents~\eqref{gg}.
A similar Zariski open set of odd elements $\mathcal W_u$ can be described as follows.
Take any \emph{odd} number of indices~$t$ and make the interchange $\partial_t\leftrightarrow x_t$ in the
equations~\eqref{dandx} \emph{for these~$t$}.
An easy exercise shows that, for a~given antisymmetric matrix~$\alpha_{tt'}^{(u)}$, the operators in the
left-hand sides of the resulting equations span an isotropic subspace, and annihilate a~one-dimensional
space of \emph{odd} $\mathcal W_u$'s.
Now take, for every~$\alpha_{tt'}^{(u)}$, those $\mathcal W_u$ that satisfy the resulting equations.

For instance, we can make the mentioned interchange for all $n=5$ indices~$t$, so that the $\mathcal W_u$'s
in the resulting set will satisfy equations
\begin{gather*}
\left(x_t-\sum_{t'\subset u}\alpha_{tt'}^{(u)}\partial_{t'}\right)\mathcal W_u=0
\qquad
\text{for all}
\quad
t\subset u.\tag*{\qed}
\end{gather*}
\renewcommand{\qed}{}
\end{proof}

\begin{theorem}\label{th:9}\qquad

\begin{enumerate}\itemsep=0pt
\item\label{th:9:i}
For weights $\mathcal W$ of the form~\eqref{gg}, both sides of~\eqref{33} satisfy $9$-dimensional spaces of
isotropic equations, i.e., equations
\begin{gather*}
d(\text{l.h.s.})=0
\qquad
\text{and}
\qquad
d(\text{r.h.s.})=0
\end{gather*}
where in both cases the relevant operators $d$ form a~$9$-dimensional isotropic space.
\item
\label{th:9:ii}
Such $9$-dimensional space of isotropic equations determines the l.h.s.\ or r.h.s.\ of~\eqref{33} up to
a~numeric factor, if it is also assumed that this l.h.s.\ or r.h.s.\ depends only on Grassmann variables on
the boundary $3$-faces, as explained after formula~\eqref{33}.
\end{enumerate}
\end{theorem}

\begin{proof}
(i)~First, we consider the integrand $\mathfrak W = \mathcal W_{12345}\mathcal W_{12346}\mathcal
W_{12356}$ for the l.h.s.\ or $\mathfrak W = \mathcal W_{12456}\mathcal W_{13456}\mathcal W_{23456}$ for
the r.h.s.
It satisf\/ies a~12-dimensional isotropic space of equations of the form
\begin{gather*}
\partial_t\mathfrak W=\sum_{t'}\gamma_{t'}x_{t'}\mathfrak W
\end{gather*}
for each boundary or inner tetrahedron $t$; these equations follow from equations~\eqref{dandx} for
individual weights and the f\/irst Leibniz rule in~\eqref{L}.
Next, if $t$ is an inner tetrahedron and if some operator $\sum\limits_{t'} (\beta_{t'}\partial_{t'} +
\gamma_{t'}x_{t'})$ anticommutes with the dif\/ferentiation $\partial_t$~-- that is,
\begin{gather}
\label{gt}
\gamma_t=0
\end{gather}
in the sum~-- and if also
\begin{gather*}
\sum_{t'}(\beta_{t'}\partial_{t'}+\gamma_{t'}x_{t'})\mathfrak W=0,
\end{gather*}
then $\partial_t \mathfrak W$ satisf\/ies a~similar equation, from which $\partial_t$ and $x_t$ are absent,
namely
\begin{gather*}
\sum_{t'\ne t}(\beta_{t'}\partial_{t'}+\gamma_{t'}x_{t'})(\partial_t\mathfrak W)=-\partial_t\sum_{t'}
(\beta_{t'}\partial_{t'}+\gamma_{t'}x_{t'})\mathfrak W=0,
\end{gather*}
and it is not hard to see that $\int \mathfrak W \,\mathrm dx_t$~-- the right derivative~-- satisf\/ies the
same equation.

Due to condition~\eqref{gt}, now there remains, at least, an 11-dimensional isotropic space of equations
instead of the 12-dimensional.
Proceeding this way further with the two remaining inner tetrahedra $t$, we get, at least, a~9-dimensional
space of equations.
As an isotropic subspace in a~18-dimensional complex Euclidean space (nine boundary tetrahedra~$t$,
operators~$\partial_t$ and~$x_t$ for each of them) cannot be more than 9-dimensional, it is exactly
9-dimensional.

(ii)~This is proved in full analogy with similar statement in item~\ref{th:5:ii} in
Theorem~\ref{th:5}.
\end{proof}

\begin{remark}
\label{r:parity}
This time, each side of~\eqref{33} is easily shown to be an \emph{odd} element~-- namely, of Grassmann
degree~3, and this determines the connected component in the manifold of maximal isotropic subspaces where
our subspaces belong.
This will be important for the construction in Section~\ref{s:p}, see Remark~\ref{r:+}.
\end{remark}

\section[A large family of Grassmann-Gaussian weights satisfying relation 3-3]{A large family of Grassmann--Gaussian weights\\ satisfying relation 3--3}
\label{s:p}

In this section, we construct a~18-parameter family of Grassmann weights depending on the variables
$x_{ijkl}$ on the boundary tetrahedra of either l.h.s.\ or r.h.s.\ of Pachner move 3--3 and such that
a~weight in this family can be represented as \emph{both} the l.h.s.\ and r.h.s.\ of~\eqref{33}, with all
the 4-simplex weights $\mathcal W_{ijklm}$ having the form~\eqref{gg}.
Although the search for an \emph{algebraic-topologically meaningful} parameterization for these weights is
still in progress, the very existence of such family is already of interest; moreover, some properties of
these weights can be seen already from the parameterization given below.

\subsection{Heuristic parameter count}

Before presenting our construction below in Subsection~\ref{ss:rig}, we would like to explain it
heuristically, using parameter counting.
For a~single 4-simplex, the corresponding isotropic space of operators, spanned by the operators in big
parentheses in~\eqref{dandx}, depends on 10 parameters.
When we compose the l.h.s.\ or r.h.s.\ of~\eqref{33} (not yet demanding that l.h.s.\ be equal to r.h.s.),
there are thus 30 parameters.
Three of them are, however, redundant, because of the possible scalings of variables $x_t$ on three
\emph{inner} tetrahedra $t$~-- it is easily seen that such scalings may only multiply the considered
integrals by a~numeric factor.
So, we have $3\times 10 - 3 = 27$ essential parameters in each side of~\eqref{33}.

On the other hand, a~9-dimensional isotropic subspace in an~18-dimensional complex Euclidean space is
determined by 36 parameters.
So, requiring this equalness, we subtract 36 parameters and are left with $2\times 27 - 36 = 18$ parameters.

\subsection{Rigorous construction}
\label{ss:rig}

There are nine boundary tetrahedra in the l.h.s.\ or r.h.s.\ of Pachner move 3--3.
We see it convenient to arrange them in the following table, where also 4-simplices are indicated by small
numbers to which the tetrahedra belong:
\begin{gather}\label{table}
\begin{array}{@{}c|ccc}
&\scriptstyle12456&\scriptstyle13456&\scriptstyle23456\\
\hline
\scriptstyle12345&1245&1345&2345\\
\scriptstyle12346&1246&1346&2346\\
\scriptstyle12356&1256&1356&2356
\end{array}
\end{gather}
Thus, the tetrahedra in every row correspond to a~4-simplex in the l.h.s., and the tetrahedra in every
column correspond to a~4-simplex in the r.h.s.\ of the move.

First, we introduce nine nonvanishing parameters $\varkappa_t\in \mathbb C$ for all tetrahedra $t=ijkl$ in
the table, and then eighteen \emph{orthonormal} vectors-operators~-- a~pair
\begin{gather*}
\mathsf e_t=\frac{1}{\varkappa_t}\frac{\partial}{\partial x_t}+\varkappa_t x_t,
\qquad
\mathsf f_t=\mathrm i\left(\frac{1}{\varkappa_t}\frac{\partial}{\partial x_t}-\varkappa_t x_t\right)
\end{gather*}
for each $t$; here
\begin{gather*}
\mathrm i=\sqrt{-1}.
\end{gather*}
Then, we introduce six more parameters: $\lambda_u$ for each table row, and $\mu_u$ for each table column,
where $u=ijklm$ is the corresponding 4-simplex.
With these parameters, we construct the following six \emph{isotropic and mutually orthogonal} vectors:
\begin{gather}
\label{g}
\mathsf g_{ijklm}=\mathsf e_{ijlm}+\mathrm i\cos\lambda_{ijklm}\,\mathsf e_{iklm}
+\mathrm i\sin\lambda_{ijklm}\,\mathsf e_{jklm}
\end{gather}
for the table rows, and
\begin{gather}
\label{h}
\mathsf h_{ijklm}=\mathsf f_{ijkl}+\mathrm i\cos\mu_{ijklm}\,\mathsf f_{ijkm}+\mathrm i\sin\mu_{ijklm}
\,\mathsf f_{ijlm}
\end{gather}
for the table columns.

Next, we bring into consideration six more unit vectors, orthogonal to each other and to all $\mathsf g_u$
and $\mathsf h_u$:
\begin{gather}
\label{p}
\mathsf p_{ijklm}=\sin\lambda_{ijklm}\,\mathsf e_{iklm}-\cos\lambda_{ijklm}\,\mathsf e_{jklm}
\end{gather}
for each row and
\begin{gather}
\label{q}
\mathsf q_{ijklm}=\sin\mu_{ijklm}\,\mathsf f_{ijkm}-\cos\mu_{ijklm}\,\mathsf f_{ijlm}
\end{gather}
for each column, and an orthogonal $3\times 3$ matrix
\begin{gather*}
A=
\begin{pmatrix}
\cos\psi&\sin\psi&0
\\
-\sin\psi&\cos\psi&0
\\
0&0&1
\end{pmatrix}
\begin{pmatrix}
1&0&0
\\
0&\cos\psi'&\sin\psi'
\\
0&-\sin\psi'&\cos\psi'
\end{pmatrix}
\begin{pmatrix}
\cos\psi''&\sin\psi''&0
\\
-\sin\psi''&\cos\psi''&0
\\
0&0&1
\end{pmatrix},
\end{gather*}
where $\psi$, $\psi'$ and $\psi''$~-- Euler angles for $A$~-- are our three remaining parameters.
With these vectors and matrix, we construct isotropic, and orthogonal to each other as well as to all~$\mathsf g_u$ and~$\mathsf h_u$, vectors~$\mathsf r$,~$\mathsf s$ and~$\mathsf t$.
It is convenient for us to arrange these vectors in a~column, and we def\/ine them as follows:
\begin{gather}
\label{rst}
\begin{pmatrix}
\mathsf r
\\
\mathsf s
\\
\mathsf t
\end{pmatrix}
=
\begin{pmatrix}
\mathsf p_{12345}
\\
\mathsf p_{12346}
\\
\mathsf p_{12356}
\end{pmatrix}
+\mathrm i A
\begin{pmatrix}
\mathsf q_{12456}
\\
\mathsf q_{13456}
\\
\mathsf q_{23456}
\end{pmatrix}.
\end{gather}

\begin{remark}\label{r:+}
The plus sign before the second term in~\eqref{rst} cannot be changed to minus without making change(s)
elsewhere in our construction.
As a~direct calculation shows, this sign ensures that the isotropic space spanned by vectors $\mathsf g_u$,
$\mathsf h_u$, $\mathsf r$, $\mathsf s$ and $\mathsf t$ (see item~\ref{th:9c:i} below in
Theorem~\ref{th:9c}) belongs to the desired connected component, according to Remark~\ref{r:parity}.
\end{remark}
\begin{theorem}\label{th:9c}\qquad
\begin{enumerate}
\itemsep=0pt
\item
\label{th:9c:i}
The linear space $\mathfrak V$ spanned by vectors $\mathsf g_{12345}$, $\mathsf g_{12346}$, $\mathsf
g_{12356}$, $\mathsf h_{12456}$, $\mathsf h_{13456}$, $\mathsf h_{23456}$, $\mathsf r$, $\mathsf s$ and~$\mathsf t$ is $9$-dimensional isotropic~-- a~maximal isotropic subspace in the $18$-dimensional space of
operators~\eqref{bg} for tetrahedra $t$ in the table~\eqref{table}.

\item
\label{th:9c:ii}
The $18$ parameters $\varkappa_t$, $\lambda_u$, $\mu_u$, $\psi$, $\psi'$ and $\psi''$, used in our
construction, are independent: the Jacobian matrix of the mapping from the space of these parameters to the
Grassmannian $($which consists of $9$-dimensional subspaces in the mentioned $18$-dimensional linear space$)$ has
rank $18$ in a~generic point.

\item
\label{th:9c:iii}
For generic parameters, $\mathfrak V$ is such that there exist such weights $\mathcal W_u$ of the
form~\eqref{gg} for all $4$-simplices in the l.h.s.\ and r.h.s.\ of~\eqref{33} that both sides of~\eqref{33}
are turned into zero by all operators in~$\mathfrak V$.
\end{enumerate}
\end{theorem}

\begin{proof}
(i)~This follows directly from our construction.

(ii)~This is shown by a~direct calculation (enough to f\/ind that rank is 18 for some
specif\/ic values of parameters).

(iii)~We begin with considering the four vectors $\mathsf g_{12346}$, $\mathsf g_{12356}$,
$\mathsf s$ and $\mathsf t$, see~\eqref{g} and~\eqref{rst}.
They are linearly independent, and their expansions in terms of the basis vectors $\mathsf e_t$ and
$\mathsf f_t$ have zero coef\/f\/icients if $t$ belongs to the f\/irst row of table~\eqref{table}; we
visualize this fact by saying that our four vectors have the form {\small $
\begin{pmatrix}
0 & 0 & 0
\\
\ast & \ast & \ast
\\
\ast & \ast & \ast
\end{pmatrix}
{}$}.
Moreover, $\mathsf g_{12346}$ and $\mathsf g_{12356}$ have the forms {\small $
\begin{pmatrix}
0 & 0 & 0
\\
\ast & \ast & \ast
\\
0 & 0 & 0
\end{pmatrix}
{}$} and {\small $
\begin{pmatrix}
0 & 0 & 0
\\
0 & 0 & 0
\\
\ast & \ast & \ast
\end{pmatrix}
{}$}; in this proof, we call vectors of such forms \emph{second-row} and \emph{third-row} vectors,
respectively.

Next, we consider the orthogonal projections of linear combinations $\sigma\mathsf s+\tau\mathsf t$ (where
$\sigma,\tau\in\mathbb C$) onto the space of third-row vectors.
There exist two such linear combinations
\begin{gather}
\label{1236dx}
\sigma_d\mathsf s+\tau_d\mathsf t
\qquad
\text{and}
\qquad
\sigma_x\mathsf s+\tau_x\mathsf t
\end{gather}
whose projections, called $\mathsf d_{1236}^{(12356)}$ and $\mathsf x_{1236}^{(12356)}$ (where the
4-simplex $12356$ corresponds to the third row, and the \emph{inner} tetrahedron $1236$ is common for it
and the ``second-row'' 4-simplex $12346$), satisfy
\begin{gather*}
\Big\langle\mathsf d_{1236}^{(12356)},\,\mathsf d_{1236}^{(12356)}\Big\rangle
=\Big\langle\mathsf x_{1236}^{(12356)},\,\mathsf x_{1236}^{(12356)}\Big\rangle=0,
\qquad
\Big\langle\mathsf d_{1236}^{(12356)},\,\mathsf x_{1236}^{(12356)}\Big\rangle=1.
\end{gather*}
As $\mathsf s$ and $\mathsf t$ lie in an isotropic subspace, this means also that the projections $\mathsf
d_{1236}^{(12346)}$ and $\mathsf x_{1236}^{(12346)}$ of the same vectors~\eqref{1236dx} onto the space of
\emph{second}-row vectors satisfy
\begin{gather*}
\Big\langle\mathsf d_{1236}^{(12346)},\,\mathsf d_{1236}^{(12346)}\Big\rangle
=\Big\langle\mathsf x_{1236}^{(12346)},\,\mathsf x_{1236}^{(12346)}\Big\rangle=0,
\qquad
\Big\langle\mathsf d_{1236}^{(12346)},\,\mathsf x_{1236}^{(12346)}\Big\rangle=-1.
\end{gather*}

Now the three operators
\begin{gather}
\label{1236(5)}
\frac{\partial}{\partial x_{1236}}+\mathsf d_{1236}^{(12356)},
\qquad
x_{1236}-\mathsf x_{1236}^{(12356)}
\qquad
\text{and}
\qquad
\mathsf g_{12356}
\end{gather}
span a~three-dimensional isotropic subspace in the (10-dimensional) space of operators~\eqref{bg} for which
$t\in \partial(12356)$ (boundary of the 4-simplex $12356$), while the three operators
\begin{gather}
\label{1236(4)}
\frac{\partial}{\partial x_{1236}}+\mathsf d_{1236}^{(12346)},
\qquad
x_{1236}+\mathsf x_{1236}^{(12346)}
\qquad
\text{and}
\qquad
\mathsf g_{12346}
\end{gather}
span a~three-dimensional isotropic subspace in the space of operators~\eqref{bg} for which $t\in
\partial(12346)$.

To move further, we note that our space $\mathfrak V$ depends on the vectors $\mathsf p_u$ and $\mathsf
q_u$ only modulo the six vectors $\mathsf g_u$ and $\mathsf h_u$: the def\/initions~\eqref{p} and~\eqref{q}
can be changed by adding any linear combinations of $\mathsf g_u$'s and $\mathsf h_u$'s to their right-hand
sides, and this does not af\/fect~$\mathfrak V$.
In particular, this means that each $\mathsf q_u$ can be changed to a~linear combination of itself and
$\mathsf h_u$~\eqref{h} (with the same $u=ijklm$) in such two ways that all the new $\mathsf q_u$'s will
f\/it into the pattern {\small $
\begin{pmatrix}
\ast & \ast & \ast
\\
\ast & \ast & \ast
\\
0 & 0 & 0
\end{pmatrix}
{}$} in the f\/irst case, and {\small $
\begin{pmatrix}
\ast & \ast & \ast
\\
0 & 0 & 0
\\
\ast & \ast & \ast
\end{pmatrix}
$} in the second case.
Then, in full analogy with what we have already done, we obtain the analogues of operators~\eqref{1236(5)}
and~\eqref{1236(4)} for the two remaining inner tetrahedra (namely, $1234$ and $1235$, respectively,
instead of $1236$) in the l.h.s.\ of the move 3--3 and their adjoining 4-simplices.

As a~result, we get a~f\/ive-dimensional isotropic space of operators for a~4-simplex in the l.h.s.
For instance, for $12356$, it is spanned by the operators~\eqref{1236(5)} and also
\begin{gather*}
\frac{\partial}{\partial x_{1235}}+\mathsf d_{1235}^{(12356)}
\qquad
\text{and}
\qquad
x_{1235}+\mathsf x_{1235}^{(12356)}
\end{gather*}
(the signs before $\mathsf x_t^{(u)}$ must be dif\/ferent for the two 4-simplices $u$ containing the
tetrahedron~$t$).

Of course, the f\/ive-dimensional isotropic space of operators can be found in a~similar way also for the
4-simplices in the r.h.s.\ of our Pachner move.
Then we def\/ine the weights $\mathcal W_u$ for all the six 4-simplices as satisfying each the
corresponding f\/ive equations, according to item~\ref{th:5:ii} in Theorem~\ref{th:5}, and (see the proof
of item~\ref{th:9:ii} in Theorem~\ref{th:9}) the relation~\eqref{33} does hold.

What remains is to show that the above construction can be performed in such way that we get \emph{even}
Grassmann elements~-- exponents~\eqref{gg} (see item~\ref{th:5:iii} in Theorem~\ref{th:5}) as our weights~$\mathcal W_u$.
We think that the easiest way to do this is to refer to the nontrivial example given below in
Subsection~\ref{ss:tan}, where the weights have the form~\eqref{gg}, and the construction of isotropic
spaces mentioned in this proof works well; then it is extended to the general case by continuity.
So, to within this example, Theorem~\ref{th:9c} is proven.
\end{proof}

\section[Two explicit constructions of parameterized weights satisfying relation 3-3]{Two explicit constructions of parameterized weights\\ satisfying relation 3--3}
\label{s:33spe}

The f\/irst construction, given in Subsection~\ref{ss:tan}, is a~particular case of the construction in the
previous Section~\ref{s:p}.
It depends on f\/ive (six, one of which is redundant) parameters~-- so, the nature of remaining $18-5=13$
parameters is still mysterious.
Note that the rank of mat\-rix~$\big(\alpha_{t_1t_2}^{(u)}\big)$ (compare formula~\eqref{gg}
with~\eqref{Phi},~\eqref{tg} and~\eqref{exp}) is 4 for this construction, as is the rank of a~generic
antisymmetric $5\times 5$ matrix, hence this rank is also 4 in the general case of Section~\ref{s:p}.

The second construction, given in Subsection~\ref{ss:mh}, is probably a~\emph{limiting} case of the
construction in Section~\ref{s:p}, because here $\rank (\alpha_{t_1t_2}^{(u)}) = 2$.
Despite this kind of degeneracy, the second construction exhibits extremely interesting relations to exotic
homologies, studied in Sections~\ref{s:h} and~\ref{s:24spe}.

\subsection{Some formulas common for the two constructions}

We are going to present, in Subsections~\ref{ss:tan} and~\ref{ss:mh}, two explicit constructions of nicely
para\-me\-terized Grassmann four-simplex weights $\mathcal W_{ijklm}$ of the form~\eqref{gg}, satisfying the
3--3 algebraic relation~\eqref{33}.
In this subsection, we write out some formulas that belong to both of them.

First, in both cases a~quantity $\varphi_{ijk}$ is introduced for each 2-face $ijk$.
These $\varphi_{ijk}$ enter both in the expressions for weights and in the multiplier $\mathrm{const}$
in~\eqref{33}.
To be more exact, our relations here look as follows:
\begin{gather}
\frac{1}{\varphi_{123}}\iiint\mathcal W_{12345}\mathcal W_{12346}\mathcal W_{12356}\,\mathrm dx_{1234}
\,\mathrm dx_{1235}\,\mathrm dx_{1236}
\nonumber
\\
\qquad
{}=-\frac{1}{\varphi_{456}}\iiint\mathcal W_{12456}\mathcal W_{13456}\mathcal W_{23456}\,\mathrm dx_{1456}
\,\mathrm dx_{2456}\,\mathrm dx_{3456}.
\label{33'}
\end{gather}

Second, the following Grassmannian quadratic form is used in both cases.
For a~4-simplex $u=ijklm$, let $abc$ be its 2-face, and $d_1<d_2$~-- two remaining vertices.
We put
\begin{gather}
\Phi_{ijklm}=\sum_{\substack{\text{ over 2-faces }abc
\\[.3ex]
\text{ of }ijklm}}\epsilon_{d_1abcd_2}\,\varphi_{abc}\,x_{\{abcd_1\}}x_{\{abcd_2\}},
\label{Phi}
\end{gather}
where $\epsilon_{d_1abcd_2}=1$ if the order $d_1abcd_2$ of vertices determines the orientation of $ijklm$
induced by the f\/ixed orientation of the manifold~-- l.h.s.\ or r.h.s.\ of a~Pachner move in our case~--
and $\epsilon_{d_1abcd_2}=-1$ otherwise.
Recall also Convention~\ref{cnv:set} concerning the curly brackets in the subscripts in~\eqref{Phi}.
\begin{remark}
In practical calculations, we use formula
\begin{gather*}
\epsilon_{d_1abcd_2}=p_{ijklm}\epsilon_{d_1abcd_2}^{ijklm},
\end{gather*}
where $p_{ijklm}$ ref\/lects the consistent orientation of 4-simplices.
Namely, for the simplices in the l.h.s.\ of move 3--3, $p_{12345}=-p_{12346}=p_{23456}=1$, and for the
simplices in the r.h.s.\ $p_{12456}=-p_{13456}=p_{23456}=1$.
As for $\epsilon_{d_1abcd_2}^{ijklm}$, it is the sign of permutation between the sequences of its
subscripts and superscripts.
\end{remark}

\subsection{First family of weights}
\label{ss:tan}

Let a~complex number $\xi_i$ be put in correspondence to every vertex $i=1,\dots,6$.
We call these numbers \emph{vertex coordinates}.
Then we def\/ine $\varphi_{ijk}$ as follows:
\begin{gather}
\label{tg}
\varphi_{ijk}=\frac{\xi_i-\xi_j}{1+\xi_i\xi_j}\cdot\frac{\xi_j-\xi_k}{1+\xi_j\xi_k}\cdot\frac{\xi_k-\xi_i}
{1+\xi_k\xi_i},
\end{gather}
then $\Phi_{ijklm}$ according to~\eqref{Phi}, and then the weight $\mathcal W_{ijklm}$ as the following
Grassmann--Gaussian exponent:
\begin{gather}
\label{exp}
\mathcal W_{ijklm}=\exp\Phi_{ijklm}.
\end{gather}

These formulas for weights f\/irst appeared in~\cite[Appendix]{2-cocycles}.
\begin{theorem}
The weights $\mathcal W_{ijklm}$ with $\varphi_{ijk}$ defined according to~\eqref{tg} satisfy the
relation~\eqref{33'}.
\end{theorem}
\begin{proof}
Direct computer calculation.
\end{proof}

Direct calculations show that the isotropic spaces of operators for the weights introduced in this
Subsection f\/it well into the scheme of Section~\ref{s:p}.
We think that this subject of isotropic spaces for 4-simplex weights deserves a~detailed study in further
works; right here we write out, just for illustration, the elegant explicit formulas for one $\mathsf g_u$
and one $\mathsf h_u$ (and, looking at them, it will not be hard to guess the formulas for other $\mathsf
g_u$ and $\mathsf h_u$).
First, introduce the following auxiliary quantities:
\begin{gather*}
r_{ijkl}
=\frac{(\xi_i\xi_j+1)(\xi_k\xi_l+1)(\xi_i\xi_j\xi_k\xi_l-\xi_k\xi_l+\xi_j\xi_l+\xi_i\xi_l+\xi_j\xi_k+\xi_i\xi_k-\xi_i\xi_j+1)}{(\xi_k-\xi_i)(\xi_k-\xi_j)(\xi_l-\xi_i)(\xi_l-\xi_j)},
\\
s_{ijkl}
=-\frac{(\xi_j-\xi_i)(\xi_l-\xi_k)(\xi_j\xi_k\xi_l+\xi_i\xi_k\xi_l-\xi_i\xi_j\xi_l-\xi_i\xi_j\xi_k+\xi_l+\xi_k-\xi_j-\xi_i)}{(\xi_i\xi_k+1)(\xi_j\xi_k+1)(\xi_i\xi_l+1)(\xi_j\xi_l+1)}.
\end{gather*}
Then, $\mathsf g_{12345}$ is proportional to the following vector:
\begin{gather}
\mathsf g_{12345} \propto r_{1245}\frac{\partial}{\partial x_{1245}}+s_{1245}x_{1245}+r_{1345}
\frac{\partial}{\partial x_{1345}}-s_{1345}x_{1345}
\nonumber
\\
\hphantom{\mathsf g_{12345} \propto}{}
+r_{2345}\frac{\partial}{\partial x_{2345}}+s_{2345}x_{2345},
\label{g1}
\end{gather}
and $\mathsf h_{12456}$ is proportional to the following vector:
\begin{gather}
\mathsf h_{12456} \propto r_{1245}\frac{\partial}{\partial x_{1245}}-s_{1245}x_{1245}+r_{1246}
\frac{\partial}{\partial x_{1246}}+s_{1246}x_{1246}
\nonumber
\\
\hphantom{\mathsf h_{12456} \propto}{}
+r_{1256}\frac{\partial}{\partial x_{1256}}-s_{1256}x_{1256}.
\label{h1}
\end{gather}

\subsection{Second family of weights}
\label{ss:mh}

This time, let each vertex $i$ have \emph{three} complex coordinates $\xi_i$, $\eta_i$, $\zeta_i$ over f\/ield~$\mathbb C$.
We def\/ine~$\varphi_{ijk}$ as the following determinant:
\begin{gather}
\label{phi}
\varphi_{ijk}=\left|
\begin{matrix}
\xi_i&\xi_j&\xi_k
\\
\eta_i&\eta_j&\eta_k
\\
\zeta_i&\zeta_j&\zeta_k
\end{matrix}
\right|.
\end{gather}

Then we def\/ine the quantity
\begin{gather}
\label{h33}
h_{ijklm}=\alpha\xi_n+\beta\eta_n+\gamma\zeta_n,
\end{gather}
where $\alpha,\beta,\gamma \in \mathbb C$, and $n \in \{1,\dots, 6\}$ is the number \emph{missing} in the set
$\{i,j,k,l,m\}$.

Finally, we def\/ine the 4-simplex weight $\mathcal W_{ijklm}$ as follows:
\begin{gather}
\label{hPhi}
\mathcal W_{ijklm}=h_{ijklm}+\Phi_{ijklm}.
\end{gather}
\begin{theorem}
The weight~\eqref{hPhi} is a~Grassmann--Gaussian exponent:
\begin{gather}
\label{he}
\mathcal W_{ijklm}=h_{ijklm}\exp(\Phi_{ijklm}/h_{ijklm}).
\end{gather}
\end{theorem}
\begin{proof}
A direct calculation shows that the form $\Phi_{ijklm}$ has now rank 2.
So, the Grassmann exponent in~\eqref{he} cannot include terms of degree${}>2$, and the terms of
degree${}\le 2$ are exactly as in~\eqref{hPhi}.
\end{proof}
\begin{theorem}
\label{th:33}
The weights $\mathcal W_{ijklm}$ defined in this Subsection satisfy the $3$--$3$ relation~\eqref{33'}.
\end{theorem}

\begin{proof}\sloppy
Direct calculation.
We used our package PL\footnote{Korepanov A.I., Korepanov~I.G., Sadykov N.M., PL: Piecewise-linear topology using GAP, \url{http://sf.net/projects/plgap/}.} for manipulations
in Grassmann al\-gebra.$\!\!\!\!\!$
\end{proof}

\begin{remark}
\label{r:r}
Formula~\eqref{h33} is simple and works well, but conceals the real nature of quanti\-ties~$h_{ijklm}$.
This will be explained below in Sections~\ref{s:h} and~\ref{s:24spe}, and we will rephrase the
statement of Theorem~\ref{th:33} in new terms, as part of Theorem~\ref{th:24}.
\end{remark}

\section{An exotic analogue of middle homologies}
\label{s:h}

The terms $h_{ijklm}$ of zero Grassmann degree in weights~\eqref{hPhi} have actually an exotic homological
nature.
This becomes especially clear if we consider not only move 3--3, but also move 2--4, and this we are going
to do in Section~\ref{s:24spe}.
Right here, we are presenting the sequence~\eqref{g3g4} of two linear mappings and some related notions and
statements.
Then, in the end of Subsection~\ref{ss:h}, we will explain that what is essential in a~set of terms
$h_{ijklm}$ corresponds to an element in $\Ker g_4 / \Image g_3$, at least in the case of the mentioned
Pachner moves.
To be exact, the permitted $h_{ijklm}$'s will correspond to $\Ker g_4$, while the result for the Grassmann
weight of either side of the move will not change under a~change of the $h_{ijklm}$'s corresponding to an
element in $\Image g_3$.

Sequence~\eqref{g3g4} is expected to be part of a~longer chain complex, but we don't need here that complex
in full.

\smallskip

Let there be an oriented triangulated PL manifold $M$ with boundary.
We introduce $\mathbb C$-linear spaces $\mathbb C^{N'_1}$ and $\mathbb C^{N_4}$ whose bases are
\emph{inner} edges and (all) 4-simplices of $M$, respectively (notations in accordance with
Convention~\ref{cnv:num}), and two $\mathbb C$-linear mappings between them as follows:
\begin{gather}
\label{g3g4}
\mathbb C^{N'_1}\stackrel{g_3}{\longrightarrow}\mathbb C^{N_4}\stackrel{g_4}{\longrightarrow}
\mathbb C^{N'_1}.
\end{gather}
We use notations $g_3$ and $g_4$ because these mappings have a~clear analogy with mappings~$g_3$ and~$g_4$
in~\cite[formula(13)]{exo}; we leave the def\/inition and discussion of other $g_i$ (namely, $g_1$, $g_2$,
$g_5$ and~$g_6$) for further papers.

A set of admissible values for $h_{ijklm}$ will correspond to an element of $\Ker g_4$, while such
Grassmann weights as the right-hand side of formula~\eqref{24} below do not change when an element of
$\Image g_3$ is added to it.

By def\/inition, the matrix element of mapping $g_3$ between an edge $b=ij$ and a~4-simplex $u$ vanishes
unless $b\subset u$.
Assuming $b\subset u$, we can write $u=\{ijklm\}$, which means, according to Convention~\ref{cnv:set}, that
$u$ has vertices $i$, $i$, $k$, $l$ and $m$ but they all don't necessarily go in the increasing order.
In this case, the matrix element of $g_3$ is
\begin{gather}
\label{g3}
\frac{1}{\varphi_{ijk}\varphi_{ijl}\varphi_{ijm}}.
\end{gather}
Recall that $\varphi_{ijk}$ is def\/ined in~\eqref{phi}.

Similarly, by def\/inition, the matrix element of mapping $g_4$ between a~4-simplex $u$ and an edge $b$ is
nonzero only if $b\subset u$.
We write again $u=\{ijklm\}$ and $b=ij$, and def\/ine this matrix element as
\begin{gather}
\label{g4}
\epsilon_{ijklm}\varphi_{klm},
\end{gather}
where $\epsilon_{ijklm}=1$ if the sequence $ijklm$, in this order, gives the consistent (with the f\/ixed
orientation of manifold $M$) orientation of $u$, and $\epsilon_{ijklm}=-1$ otherwise.
\begin{theorem}
\label{th:c}
Mappings $g_3$ and $g_4$ defined according to~\eqref{g3} and~\eqref{g4} form a~chain, in the sense that
\begin{gather*}
g_4\circ g_3=0.
\end{gather*}
\end{theorem}

\begin{proof}
We f\/irst show this for $M=\partial\Delta^5$~-- the sphere $S^4$ triangulated into f\/ive 4-simplices as
the boundary of 5-simplex $123456$.
Consider an edge $a\subset S^4$ as a~basis vector in the leftmost space in~\eqref{g3g4}, an edge $b\subset
S^4$ as a~basis vector in the rightmost space in~\eqref{g3g4}, and the matrix element of mapping $g_4\circ
g_3$ between $a$ and $b$.
This matrix element is the sum of products of expressions~\eqref{g3} and~\eqref{g4} over those 4-simplices
$u$ that contain both $a$ and $b$; we call such a~product \emph{contribution} of $u$.

There are three cases; in all three $ijklmn$ is a~permutation of the set of vertices $1,\dots, 6$.

Case~1. Edges $a$ and $b$ coincide, $a=b=ij$.
The matrix element is
\begin{gather*}
\frac{\epsilon_{ijklm}\varphi_{klm}}{\varphi_{ijk}\varphi_{ijl}\varphi_{ijm}}+\frac{\epsilon_{ijkln}
\varphi_{kln}}{\varphi_{ijk}\varphi_{ijl}\varphi_{ijn}}+\frac{\epsilon_{ijkmn}\varphi_{kmn}}{\varphi_{ijk}
\varphi_{ijm}\varphi_{ijn}}+\frac{\epsilon_{ijlmn}\varphi_{lmn}}{\varphi_{ijl}\varphi_{ijm}\varphi_{ijn}}
\\
\qquad
{}=\epsilon_{ijklm}\frac{\varphi_{klm}\varphi_{ijn}-\varphi_{kln}\varphi_{ijm}+\varphi_{kmn}\varphi_{ijl}
-\varphi_{lmn}\varphi_{ijk}}{\varphi_{ijk}\varphi_{ijl}\varphi_{ijm}\varphi_{ijn}}=0.
\end{gather*}
The numerator, after the reducing to a~common denominator, vanishes: this is a~\emph{Pl\"ucker relation}.

Case~2. Edges $a$ and $b$ have one vertex in common, $a=ij$, \;$b=ik$.
The matrix element is
\begin{gather*}
\frac{\epsilon_{ikjlm}\varphi_{jlm}}{\varphi_{ijk}\varphi_{ijl}\varphi_{ijm}}+\frac{\epsilon_{ikjln}
\varphi_{jln}}{\varphi_{ijk}\varphi_{ijl}\varphi_{ijn}}+\frac{\epsilon_{ikjmn}\varphi_{jmn}}{\varphi_{ijk}
\varphi_{ijm}\varphi_{ijn}}
=\epsilon_{ikjlm}\frac{\varphi_{jlm}\varphi_{ijn}-\varphi_{jln}\varphi_{ijm}+\varphi_{jmn}\varphi_{ijl}}
{\varphi_{ijk}\varphi_{ijl}\varphi_{ijm}\varphi_{ijn}}=0.
\end{gather*}
Here the numerator is obtained from the numerator in our Case 1 by letting $k=j$.

Case~3. Edges $a$ and $b$ do not intersect, $a=ij$, \;$b=kl$.
The matrix element is
\begin{gather*}
\frac{\epsilon_{klijm}\varphi_{ijm}}{\varphi_{ijk}\varphi_{ijl}\varphi_{ijm}}+\frac{\epsilon_{klijn}
\varphi_{ijn}}{\varphi_{ijk}\varphi_{ijl}\varphi_{ijn}}=0.
\end{gather*}

Next, we prove the theorem again for $M=S^4$, but triangulated in an \emph{arbitrary} way.
This arbitrary triangulation can be achieved by a~sequence of Pachner moves performed on the initial
triangulation considered above.
Each Pachner move replaces some 4-simplices with some other ones, in such way that the withdrawn and the
replacing 4-simplices will form together a~sphere~$\partial \Delta^5$ \emph{if we change the orientation}
of, say, the withdrawn 4-simplices.
It follows then from what we have proved for~$\partial \Delta^5$ that the contribution of all the replacing
4-simplices into any matrix element of~$g_4\circ g_3$ is the same as the contribution of all the withdrawn
4-simplices, \emph{including the cases where an edge appears or disappears} during the move.

Finally, it remains to say that any manifold $M$ has locally the same structure as~$S^4$, and any matrix
element of $g_4\circ g_3$ consists only of obviously local contributions.
\end{proof}

\begin{xpr}
\label{x:r}
For a~closed oriented 4-dimensional PL manifold $M$, the vector space $\Ker g_4 / \Image g_3$ is six times
$($i.e., isomorphic to the direct sum of six copies of$)$ usual second homologies $H_2(M;\mathbb C)$.
\end{xpr}

 This has been checked (in particular) for $M=T^4$, $T^2\times S^2$, $S^1\times S^3$, $S^2\times
S^2$, $S^4$, and the Kummer surface.

\section[Relation 2-4 and the exotic homological nature of $h_{ijklm}$]{Relation 2--4 and the exotic homological nature of $\boldsymbol{h_{ijklm}}$}  \label{s:24spe}

As we have already said, it is the relation 2--4 that makes especially clear the exotic homological nature
of free terms $h_{ijklm}$.
Here it is:
\begin{gather}
\int\mathcal W_{12345}\mathcal W_{12346}\,\mathrm dx_{1234}
=-\frac{1}{\varphi_{156}\varphi_{256}\varphi_{356}\varphi_{456}}\idotsint\mathcal W_{12356}
\mathcal W_{12456}\mathcal W_{13456}\mathcal W_{23456}
\nonumber
\\
\hphantom{\int\mathcal W_{12345}\mathcal W_{12346}\,\mathrm dx_{1234}=}{}
\times w_{56}\,\mathrm dx_{1256}\,\mathrm dx_{1356}\,\mathrm dx_{1456}\,\mathrm dx_{2356}\,\mathrm dx_{2456}
\,\mathrm dx_{3456}.
\label{24}
\end{gather}
It involves a~new factor $w_{56}$ in the integrand in its r.h.s.~-- the \emph{edge weight} of the inner edge~$56$.
So, we f\/irst def\/ine this weight in Subsection~\ref{ss:w}, and then return to our quantities~$h_{ijklm}$
in Subsection~\ref{ss:h}.

\subsection{The edge weight}\label{ss:w}

For any edge $a=ij$ in a~triangulation of a~4-manifold with boundary, we def\/ine a~Grassmann
dif\/ferential operator $\partial_a=\partial_{ij}$ as the following sum over all tetrahedra $t=\{ijkl\}$
containing this edge:
\begin{gather}
\label{dij}
\partial_{ij}=\sum_{t=\{ijkl\}}\frac{1}{\varphi_{ijk}\varphi_{ijl}}\,\partial_t .
\end{gather}
\begin{lemma-nonnumbered}
The weight $\mathcal W_u$~\eqref{hPhi} of a~$4$-simplex $u$ satisfies ``edge equations''
\begin{gather*}
\partial_a\mathcal W_u=0
\end{gather*}
for any edge $a$.
\end{lemma-nonnumbered}

\begin{proof}
Direct calculation (which is, of course, nontrivial only if $a\subset u$).
\end{proof}

Specif\/ically, here is how the operator~\eqref{dij} looks for the inner edge 56 of the cluster of four
4-simplices corresponding to the r.h.s.\ of relation~\eqref{24}:
\begin{gather}
\partial_{56}=\frac{1}{\varphi_{156}\varphi_{256}}\partial_{1256}+\frac{1}{\varphi_{156}\varphi_{356}}
\partial_{1356}+\frac{1}{\varphi_{156}\varphi_{456}}\partial_{1456}
\nonumber
\\
\phantom{\partial_{56}=}
+\frac{1}{\varphi_{256}\varphi_{356}}\partial_{2356}+\frac{1}{\varphi_{256}\varphi_{456}}\partial_{2456}
+\frac{1}{\varphi_{356}\varphi_{456}}\partial_{3456}.
\label{d56}
\end{gather}
\begin{theorem}
\label{th:w}
The right-hand side of~\eqref{24} does not change if an expression is added to $w_{56}$ whose ``edge~$56$
derivative'' is zero, as follows:
\begin{gather*}
w_{56}\mapsto w_{56}+\tilde w_{56},
\qquad
\partial_{56}\tilde w_{56}=0.
\end{gather*}
\end{theorem}

\begin{proof}
We have to prove that
\begin{gather}
\idotsint\mathcal W_{12356}\mathcal W_{12456}\mathcal W_{13456}\mathcal W_{23456}
\tilde w_{56}\,\mathrm dx_{1256}\,\mathrm dx_{1356}\,\mathrm dx_{1456}\,\mathrm dx_{2356}
\,\mathrm dx_{2456}\,\mathrm dx_{3456}=0.
\label{tw56}
\end{gather}
As the ``edge 56 derivative'' of every factor in the integrand
\begin{gather*}
\mathcal X=\mathcal W_{12356}\mathcal W_{12456}\mathcal W_{13456}\mathcal W_{23456}\tilde w_{56}.
\end{gather*}
vanishes, it vanishes for all the product:
\begin{gather}
\label{l}
\partial_{56}\mathcal X=0.
\end{gather}
It is an easy exercise (just break $\mathcal X$ into the even and odd parts) to see that the right analogue
of~\eqref{l} also holds:
\begin{gather*}
\mathcal X\overleftarrow\partial_{\!56}=0,
\end{gather*}
where $\overleftarrow \partial_{\!56}$ is the same linear combination~\eqref{d56}, but with left
derivatives replaced with the right ones.
Finally, recall that the integration means, according to Section~\ref{s:G}, the right
dif\/ferentiation, and specif\/ically in~\eqref{tw56} this can be represented as follows
(compare~\eqref{dmB}):
\begin{gather*}
\mathcal X\overleftarrow\partial_{\!1256}\overleftarrow\partial_{\!1356}\cdots\overleftarrow\partial_{\!3456}
=\varphi_{156}\varphi_{256}\mathcal X\overleftarrow\partial_{\!56}\overleftarrow\partial_{\!1356}
\cdots\overleftarrow\partial_{\!3456}=0.\tag*{\qed}
\end{gather*}
\renewcommand{\qed}{}
\end{proof}

Due to Theorem~\ref{th:w}, the following def\/inition of $w_{56}$ is not surprising:
\begin{gather}
\label{w56}
w_{56}=\partial_{56}^{-1}1,
\end{gather}
by which we understand \emph{any} Grassmann algebra element $w_{56}$ such that $\partial_{56}w_{56}=1$.
For instance, we can take $w_{56} = \varphi_{156}\varphi_{256}x_{1256}$.

\subsection[The quantities $h_{ijklm}$]{The quantities $\boldsymbol{h_{ijklm}}$}
\label{ss:h}

The quantities $h_{ijklm}$ for \emph{both} moves 3--3 and 2--4 considered in this paper can be obtained as
follows.
First, glue together both sides of a~Pachner move in the natural way~-- identifying like-named boundary
simplices.
This gives $S^4=\partial \Delta^5$ (recall that this means a~4-sphere triangulated as the boundary of
a~5-simplex).
\begin{remark}
This gluing implies that we have changed the orientation in one of the sides.
Nevertheless, our mapping $g_3$ is def\/ined in such way that the orientation issues do not af\/fect the
def\/inition of allowable set of values for $h_{ijklm}$ given below.
\end{remark}
\begin{theorem}
The sequence~\eqref{g3g4} for $S^4=\partial \Delta^5$ is exact: $\Ker g_4=\Image g_3$.
\end{theorem}
\begin{proof}
Taking into account Theorem~\ref{th:c}, it is enough to show that
\begin{gather*}
\rank g_3+\rank g_4=6,
\end{gather*}
where 6 is the number of 4-simplices in the triangulation and thus the dimension of the middle space
in~\eqref{g3g4}.
This is done by a~direct calculation; both ranks prove to be 3.
\end{proof}
\begin{remark}
Compare this also with the Experimental result on page \pageref{x:r}.
\end{remark}

We now take an arbitrary chain $c_{\rm edges}$ on edges of our $S^4$~-- element of the f\/irst linear space
in~\eqref{g3g4}, and then its image under $g_3$~-- a~chain on 4-simplices.
We call the resulting coef\/f\/icients at 4-simplices $u=ijklm$, both in the l.h.s.\ and r.h.s.\ of our
Pachner move, an \emph{allowable set} of values for $h_{ijklm}$.

\begin{theorem}\label{th:24}
Let there be an allowable set of values $h_u=h_{ijklm}$ for the six $4$-simplices in both sides of either
Pachner move $3$--$3$ or $2$--$4$, and, in case of move $2$--$4$, let there be chosen any edge weight~$w_{56}$
according to~\eqref{w56}.
Let also the $4$-simplex weights be defined according to~\eqref{hPhi}, with the quadratic forms $\Phi_u$
defined according to~\eqref{Phi} and~\eqref{phi}.
Then, the relation~\eqref{33'} holds for move $3$--$3$ or, respectively,~\eqref{24} holds for move $2$--$4$.
\end{theorem}

\begin{proof}
Direct calculation.
\end{proof}

\begin{remark}
It is an easy exercise to show that our allowable sets of values $h_u=h_{ijklm}$ for Pachner moves can be
represented in the form~\eqref{h33} (although this form may disguise their nature).
Hence, the part of Theorem~\ref{th:24} dealing with move 3--3 says the same as Theorem~\ref{th:33}, as was
promised in Remark~\ref{r:r}.
\end{remark}

We explain now why we said in the f\/irst paragraph of Section~\ref{s:h} that what is essential in
a~set of terms $h_{ijklm}$ corresponds to an element in $\Ker g_4 / \Image g_3$.
First, we note that in both sides of the move 3--3, as well as in the l.h.s.\ of move 2--4, there are
\emph{no} inner edges.
So, the sequence~\eqref{g3g4} assumes the form
\begin{gather*}
0\stackrel{g_3}{\longrightarrow}\mathbb C^{N_4}\stackrel{g_4}{\longrightarrow}0,
\end{gather*}
i.e., the factor $\Ker g_4 / \Image g_3$ obviously gives \emph{all} possible values for $h_{ijklm}$, which
agrees with Theorem~\ref{th:24}.
\begin{remark}
Of course, Theorem~\ref{th:24} tells us more, namely, also how to make the $h_{ijklm}$'s agree \emph{for
the two sides}, so that~\eqref{33'} or~\eqref{24} holds.
\end{remark}

Now consider the r.h.s.\ of the 2--4 relation~\eqref{24}.
As it is equal to the l.h.s., it does not depend on the coef\/f\/icient in $c_{\rm edges}$ at the edge
$56$, the latter being absent from the l.h.s.\ of the Pachner move, while being the only inner edge for its
r.h.s.
Consider the sequence~\eqref{g3g4} for the r.h.s.\ of move 2--4.
The f\/irst and third spaces in this sequence are one-dimensional, the only basis vector being the edge
$56$.
Any allowable set of $h_u$, with the $u$'s in this r.h.s., obviously makes a~\emph{cycle} in the sense that
it is annihilated by the mapping $g_4$.
So, the r.h.s.\ of~\eqref{24} is determined by this cycle modulo a~\emph{boundary}~-- image of $g_3$.

Especially interesting question for future research is to uncover the analogue(s) of such
exotic-homo\-lo\-gical structures for more general 4-simplex weights described in Section~\ref{s:p}.

\subsection*{Acknowledgements}

We thank the creators and maintainers of GAP\footnote{GAP~-- Groups, Algorithms, Programming~-- a~System
for Computational Discrete Algebra, \url{http://gap-system.org/}.} and Maxima\footnote{Maxima, a~computer
algebra system, \url{http://maxima.sourceforge.net}.} for their excellent computer algebra systems.
We also thank the referees for valuable comments, including drawing our attention to the literature on
holonomic functions (see Remark~\ref{r:PWZ}).

\pdfbookmark[1]{References}{ref}
\LastPageEnding


\begin{thebibliography}{99}
\footnotesize
\itemsep=0pt

\bibitem{B}
Berezin F.A., The method of second quantization, \textit{Pure and Applied
  Physics}, Vol.~24, Academic Press, New York, 1966.

\bibitem{B-super}
Berezin F.A., Introduction to superanalysis, \textit{Mathematical Physics and
  Applied Mathematics}, Vol.~9, D.~Reidel Publishing Co., Dordrecht, 1987.

\bibitem{Bjork}
Bj{\"o}rk J.E., Rings of dif\/ferential operators, \textit{North-Holland
  Mathematical Library}, Vol.~21, North-Holland Publishing Co., Amsterdam,
  1979.

\bibitem{Cartier}
Cartier P., D\'emonstration ``automatique'' d'identit\'es et fonctions
  hyperg\'eom\'etriques (d'apr\`es {D}.~{Z}eilberger), S{\'e}minaire Bourbaki,
  Vol.~1991/92, \textit{Ast\'erisque}
  \textbf{206} (1992), Exp.\ No.~746, 3, 41--91.

\bibitem{C}
Chevalley C., The algebraic theory of spinors and {C}lif\/ford algebras,
  Collected works, Vol.~2, Springer-Verlag, Berlin, 1997.


\bibitem{0806}
Korepanov I.G., Geometric torsions and an {A}tiyah-style topological f\/ield
  theory, \href{http://dx.doi.org/10.1007/s11232-009-0028-0}{\textit{Theoret. and Math. Phys.}} \textbf{158} (2009), 344--354,
  \href{http://arxiv.org/abs/0806.2514}{arXiv:0806.2514}.

\bibitem{11S}
Korepanov I.G., Relations in {G}rassmann algebra corresponding to three- and
  four-dimensional {P}achner moves, \href{http://dx.doi.org/10.3842/SIGMA.2011.117}{\textit{SIGMA}} \textbf{7} (2011), 117,
  23~pages, \href{http://arxiv.org/abs/1105.0782}{arXiv:1105.0782}.

\bibitem{exo}
Korepanov I.G., Deformation of a $3\to 3$ {P}achner move relation capturing
  exotic second homologies, \href{http://arxiv.org/abs/1201.4762}{arXiv:1201.4762}.

\bibitem{2-cocycles}
Korepanov I.G., Special 2-cocycles and 3--3 {P}achner move relations in
  {G}rassmann algebra, \href{http://arxiv.org/abs/1301.5581}{arXiv:1301.5581}.

\bibitem{KS}
Korepanov I.G., Sadykov N.M., Pentagon relations in direct sums and {G}rassmann
  algebras, \href{http://dx.doi.org/10.3842/SIGMA.2013.030}{\textit{SIGMA}} \textbf{9} (2013), 030, 16~pages,
  \href{http://arxiv.org/abs/1212.4462}{arXiv:1212.4462}.

\bibitem{Lickorish}
Lickorish W.B.R., Simplicial moves on complexes and manifolds, in Proceedings
  of the {K}irbyfest ({B}erkeley, {CA}, 1998), \href{http://dx.doi.org/10.2140/gtm.1999.2.299}{\textit{Geom. Topol. Monogr.}},
  Vol.~2, Geom. Topol. Publ., Coventry, 1999, 299--320,
  \href{http://arxiv.org/abs/math.GT/9911256}{math.GT/9911256}.

\bibitem{Pachner}
Pachner U., P.{L}.\ homeomorphic manifolds are equivalent by elementary
  shellings, \textit{European~J. Combin.} \textbf{12} (1991), 129--145.

\bibitem{PWZ}
Petkov{\v{s}}ek M., Wilf H.S., Zeilberger D., {$A=B$}, A~K~Peters Ltd.,
  Wellesley, MA, 1996.

\end{thebibliography}
\end{document}